\documentclass[11pt,reqno]{amsart}
\usepackage{amsmath,amssymb,amsthm}
\usepackage{graphicx}
\usepackage[all]{xypic}
\usepackage[utf8,nocaptions]{vietnam}
\usepackage{calrsfs}
\input xypic
\def\thmsection{section}
\def\thmchangesection{changesection}

\def\thmchangechapter{changechapter}
\def\thmchange{change}
\def\thmplain{plain}
\ifx\theoremnumberstyle\thmplain
  \theoremstyle{break-italic}
  \newtheorem{satz}{Satz}
\else
  \ifx\theoremnumberstyle\thmsection
    \theoremstyle{break-italic}
    \newtheorem{satz}{Satz}[section]
  \else
    \ifx\theoremnumberstyle\thmchange
      \swapnumbers
      \theoremstyle{break-italic}
      \newtheorem{satz}{Satz}
    \else
      \ifx\theoremnumberstyle\thmchangesection
         \swapnumbers
         \theoremstyle{break-italic}
         \newtheorem{satz}{Satz}[section]
      \else
        \ifx\theoremnumberstyle\thmchangechapter
           \swapnumbers
           \theoremstyle{break-italic}
           \newtheorem{satz}{Satz}[chapter]
        \else
          \ifx\theoremnumberstyle\thmsection
             \theoremstyle{break-italic}
             \newtheorem{satz}{Satz}[section]
          \else
            \theoremstyle{break-italic}
            \newtheorem{satz}{Satz}[section]
          \fi
        \fi
      \fi
    \fi
  \fi
\fi

\theoremstyle{break-italic}
\newtheorem{theorem}[satz]{Theorem}
\newtheorem{lemma}[satz]{Lemma}

\newtheorem{corollary}[satz]{Corollary}
\newtheorem{Proposition}[satz]{Proposition}

\newtheorem*{conjecture*}{Conjecture}

\theoremstyle{break-roman}
\newtheorem{definition}[satz]{Definition}

\newtheorem{example}[satz]{Example}

\newtheorem{remark}[satz]{Remark}

\newtheorem{conjecture}[satz]{Conjecture}

\theoremstyle{standard}

\newtheorem*{claim}{Claim}

\theoremstyle{varthm-roman}
\newtheorem*{varthm-roman}{}
\theoremstyle{varthm-italic}
\newtheorem*{varthm-italic}{}
\theoremstyle{varthm-roman-break}
\newtheorem*{varthm-roman-break}{}
\theoremstyle{varthm-italic-break}
\newtheorem*{varthm-italic-break}{}
\theoremstyle{varthm-roman-no-punctuation}
\newtheorem{varthm-roman-no-punctuation-numbered}[satz]{}
\theoremstyle{varthm-italic-no-punctuation}
\newtheorem{varthm-italic-no-punctuation-numbered}[satz]{}

\newenvironment{varthm-roman-numbered}[1]{
  \begin{varthm-roman-no-punctuation-numbered}
    \mbox{\rm\textbf{#1}}
  }{\end{varthm-roman-no-punctuation-numbered}}
\newenvironment{varthm-italic-numbered}[1]{
  \begin{varthm-italic-no-punctuation-numbered}
    \mbox{\rm\textbf{#1}}
  }{\end{varthm-italic-no-punctuation-numbered}}
\newenvironment{varthm-roman-break-numbered}[1]{
  \begin{varthm-roman-no-punctuation-numbered}
    \mbox{\rm\textbf{#1}\newline}
  }{\end{varthm-roman-no-punctuation-numbered}}
\newenvironment{varthm-italic-break-numbered}[1]{
  \begin{varthm-italic-no-punctuation-numbered}
    \mbox{\rm\textbf{#1}}\newline
  }{\end{varthm-italic-no-punctuation-numbered}}
\numberwithin{equation}{section}

\def\ex{\begin{example}
  }
  \def\eex{\end{example}}
\def\thr{\begin{theorem}}
\def\ethr{\end{theorem}}
\def\pro{\begin{Proposition}}
\def\epro{\end{Proposition}}
\def\coro{\begin{corollary}}
\def\ecoro{\end{corollary}}
\def\df{\begin{definition}}
\def\edf{\end{definition}}
\def\lm{\begin{lemma}}
\def\elm{\end{lemma}}
\def\pf{\begin{proof}}
\def\epf{\end{proof}}
\def\problem{\begin{problem}}
\def\eproblem{\end{problem}}

\def\ite{\begin{itemize}}
\def\hite{\end{itemize}}
\def\rem{\begin{remark}}
\def\erem{\end{remark}}

\def\cla{\begin{claim}}
\def\ecla{\end{claim}}
\def\conj{\begin{conjecture}}
\def\econj{\end{conjecture}}

\newcommand{\vc}{\infty}
\newcommand{\vm}{\forall}
\newcommand{\tru}{\setminus}

\newcommand{\ps}{\dfrac}
\newcommand{\can}{\sqrt}
\newcommand{\mtn}{\rightarrow}

\newcommand{\seq}[1]{\left<#1\right>}

\def\r{\Bbb R}

\def\n{\Bbb N}

\def\c{\Bbb C}

\def\z{\Bbb Z}

\begin{document}

\title[Morse polynomial functions and polynomial optimization]{Nonnegative Morse polynomial functions and polynomial optimization}
\author{L\^{e} C\^{o}ng-Tr\`{i}nh }
\address{L\^{e} C\^{o}ng-Tr\`{i}nh \\
Department of Mathematics, Quy Nhon University\\
170 An Duong Vuong, Quy Nhon, Binh Dinh}
\email{lecongtrinh@qnu.edu.vn}
\subjclass[2010]{11E25, 13J30, 14H99, 14P05, 14P10,    90C22}


\dedicatory{Dedicated to Professor H\`{a} Huy Vui  on the occasion of his 65th birthday}

\keywords{Sum of squares; Positivstellensatz;  Polynomial optimization; Local-global principle; Morse function; Non-degenerate polynomial map}

\begin{abstract} In this paper we study the  representation of Morse polynomial functions which are  nonnegative on a compact basic closed semi-algebraic set in $\r^{n}$, and having only finitely many zeros in this  set.  Following C. Bivi\`{a}-Ausina (Math Z 257:745–767, 2007), we introduce  two  classes of  non-degenerate  polynomials for which the algebraic sets defined by them are compact.  As a consequence, we study the  representation of nonnegative Morse   polynomials on  these kinds of non-degenerate algebraic sets. Moreover, we apply these results to study the polynomial optimization problem for Morse polynomial functions.
\end{abstract}

\maketitle
\tableofcontents
\section{Introduction}
Let us denote by $\r[X]$ the ring of real polynomials in $n$ variables $x_1,\cdots,x_n$, and by $\sum \r[X]^{2}$ the set of all finitely many sums of squares (SOS) of polynomials in $\r[X]$.  Let us fix a  finite subset $G=\{g_1,\cdots, g_m\}$ in   $ \r[X]$. Let   
$$K_G=\{x=(x_1,\cdots,x_n)\in \r^n| g_1(x)\geq 0,\cdots,g_m(x)\geq 0\}$$
be the \textit{basic  closed  semi-algebraic set} in $\r^n$ generated by $G$.  Let 
$$ M_G:=\{\sum_{i=1}^m s_i g_i | s_i \in \sum \r[{X}]^2\}$$
be  the \textit{quadratic module} in $\r[X]$ generated by $G$, and  let 
$$ T_G:=\{\sum_{\sigma=(\sigma_1,\cdots,\sigma_m)\in \{0,1\}^m} s_\sigma g_1^{\sigma_1}\cdots g_m^{\sigma_m}| s_\sigma \in \sum \r[X]^2\}$$ denote the \textit{preordering} in  $\r[{X}]$ generated by $G$.   It is clear that $M_G\subseteq T_G$, 
and   if a polynomial   belongs to  $T_G$ (or $M_G$)  then it is nonnegative on $K_G$. However the converse is not always true, that means there  exists a polynomial which is nonnegative on $K_G$ but it does not belong to $T_G$ (resp. $M_G$).  The well-known examples (cf. \cite{Mar08}) are Motzkin's polynomial, Robinson's polynomial, etc. in the case $G=\emptyset$ (then $K_G=\r^{n}$ and $T_G=M_G=\sum \r[X]^{2}$). 

In 1991, Schm\"{u}dgen \cite{Sch} showed, that if a polynomial is positive on a compact basic closed semi-algebraic set then it belongs to the corresponding preordering. After that, Putinar (\cite{Pu}, 1993) showed that if a  polynomial  is positive on a basic closed 
semi-algebraic set  whose associated quadratic module is Archimedean, then it  belongs to that  quadratic module.  

If we allow the polynomial $f$ having zeros in $K_G$, then the results above of Schm\"{u}dgen and  Putinar  are not true.  
Indeed, let us consider  the following counter-example, which was given by G. Stengle  \cite{Ste96}. Consider the set  $G=\{(1-x^2)^3\}$ in the ring  $\r[x]$ of  real polynomials in one variable.  For this set,   $K_G$ is the closed interval  $[-1,1]\subseteq \r$  which is compact.  The polynomial  $f=1-x^2 \in \r[x]$ is nonnegative on $[-1,1]$, and it has two zeros in $[-1,1]$. It is not difficult to show that  
$$f \notin T_G=M_G=\{s_0+s_1(1-x^2)^3 | s_0, s_1 \in \sum \mathbb{R}[x]^2\}.$$  
Therefore a natural question is that \textit{under which conditions  a polynomial which is nonnegative on a basic closed  semi-algebraic set belonging to the corresponding preordering (or quadratic module)}?  C. Scheiderer (2003 and 2005) has given the following \textit{local-global principles} to answer this question.

\thr[{\cite[Corollary 3.17]{Schei03}}]\label{local-global1} Let  $G, K_G$ and $T_G$ be as above, and let $f\in \r[X]$. Assume  that the following conditions hold true:
\ite 
\item[(1)] $K_G$ is compact;
\item[(2)]  $f\geq 0$ on $K_G$, and  $f$ has only finitely many zeros $p_1,\cdots, p_r$ in $K_G$; 
\item[(3)]  at each $p_i$, $f\in \widehat{T}_{p_i}$.
\hite
 Then $f \in T_G$.
\ethr 
Here $\widehat{T}_{p}$ (resp. $\widehat{M}_{p}$) denotes the preordering  (resp. quadratic module) generated by $T_G$  (resp. $M_G$) in the completion $\r[[X-p]]$ of the polynomial ring $\r[X]$ at the point $p\in \r^{n}$.

\thr[{\cite[Proposition 3.4]{Schei05}}]\label{local-global2} Let  $G, K_G$ and $M_G$ be as above, and let $f\in \r[X]$. Assume  that 
\ite 
\item[(1)] $M_G$ is Archimedean;
\item[(2)]  $f\geq 0$ on $K_G$, and  $f$ has only finitely many zeros $p_1,\cdots, p_r$ in $K_G$; 
\item[(3)]  at each $p_i$, $f\in \widehat{M}_{p_i}$, 
\hite
and at least one of the following conditions is satisfied:
\ite 
\item[(4)] $\dim \mathcal{V}(f) \leq 1$;
\item[(4')] for every  $p_i$, there exists a neighborhood $U$ of $p_i$ in $\r^{n}$ and an element $a\in M_G$ such that $\{a\geq 0\}\cap \mathcal{V}(f) \cap U \subseteq K_G$.
\hite 
 Then $f \in M_G$.
\ethr 

\hspace{-0.5cm} Here $\mathcal{V}(f)=\{x\in \r^{n}| f(x)=0\}$ denotes the vanishing set of $f$ in $\r^{n}$. 

The assumption on the compactness of the basic closed semi-algebraic set $K_G$ or on the Archimedean property of $M_G$ is 
necessary, and it is not difficult to verify (for Archimedean property of $M_G$, we can use, for example,  \textit{Putinar's criterion} \cite{Pu}).  However, it is complicated  and hence not convenient  in practice to verify that $f\in \widehat{T}_{p_i}$ (resp. $f\in \widehat{M}_{p_i}$) at each zero $p_i$ of $f$ in $K_G$.  Therefore, it is necessary to give a \textit{generic class} of polynomials which satisfies these conditions. 

A smooth function $f: M \mtn \r$ on a smooth manifold $M$ of dimension $n$  is called a \textit{Morse function} if all of its critical points are non-degenerate, i.e. if $p\in M$ is a critical point of $f$ then the Hessian matrix $D^{2}f(p)$ of $f$ at $p$ is invertible. 

It is well-known from Differential Topology and Singularity theory that almost smooth functions on smooth manifolds are Morse (cf. \cite{BH04}).   Furthermore,  in    Theorem \ref{thr-PosMorse} of section 2,  we show that Morse polynomial functions  solve  the disadvantage mentioned above.

The assumption on the compactness of the basic closed semi-algebraic $K_G$ in the theorems of Schm\"{u}dgen, Putinar and Scheiderer cannot be removed. In section 3 of this paper, following C. Bivi\`{a}-Ausina \cite{Bi07}, we introduce two  classes of  non-degenerate  polynomials for which the algebraic sets defined by them are compact.  As a consequence, we give a representation of nonnegative Morse   polynomials on the  non-degenerate algebraic set $K_G$ (see Corollary \ref{coro-Morsecompact1} and Corollary  \ref{coro-Morsecompact2}).

 In section 4 we give some applications of the representation of Morse polynomial functions on compact basic closed semi-algebraic sets. 
For  the global  polynomial optimization problem $$f^{*} = \min_{x\in \r^{n}} f(x),$$
 J.-B. Lasserre \cite{La01} and some other authors  have given an SOS relaxation    for this problem, which can be translated into an SDP. The finite convergence of the SOS relaxation depends  mainly on the SOS representation of $f-f^{*}$ modulo the gradient ideal $I_{grad}(f)$ of $f$. One of the sufficient conditions for the finite convergence of the SOS relaxation is that the  gradient ideal $I_{grad}(f)$ is  radical. We show in Proposition \ref{pro-sosgrad1} that for Morse polynomial functions we don't need this condition. We apply this result to show in Theorem  \ref{thr-sosconvergence1} that for Morse polynomial functions, the above SOS relaxation has a finite convergence. \\
For the constrained polynomial optimization problem on the basic closed semi-algebraic set $K_G$ 
$$ f^*= \min_{x\in K_G} f(x),$$
one of the sufficient conditions for the finite convergence of the SOS relaxation is that the KKT ideal associated to the KKT system is of dimension zero (i.e. the corresponding complex KKT variety has only finitely many points) and radical.  Then $f-f^{*}$ is  in the KKT quadratic module (resp. KKT preordering). For Morse polynomial functions, we show in Proposition  \ref{pro-kkt-representation}  that if $K_G$ is compact (resp. $M_G$ is Archimedean), and  if $f$ has only finitely many real KKT points in the interior of $K_G$, then $f-f^{*}$ belongs to $T_G$ (resp. $M_G$).  \\
J.-B. Lasserre \cite{La01} constructed a convex LMI problem in terms of the moment matrices to give a way to compute $f^{*}$ in the case where the quadratic module $M_G$ is assumed to be Archimedean.  In his method,  he assumed that  $f-f^{*}\in M_G$. Applying Proposition  \ref{pro-kkt-representation} we  can omit this assumption (see Corollary  \ref{coro-momentrepresentation}).

\hspace{-0.5cm} \textbf{Notation:}   Throughout  this paper, we denote $\r_{+}$ for  the set of nonnegative real  numbers; $\z_{+}$   the set of nonnegative integers;  $\r[[X]]$   the ring of formal power series in $n$ variables $x_1,\cdots,x_n$;  $\sum \r[X]^{2}$ (resp. $\sum \r[[X]]^{2}$)  the set of all sums of squares (SOS) of finitely many polynomials (resp. formal power series)  in $\r[X]$ (resp.   $ \r[[X]]$); $\r[[X-p]]$ the ring of formal power series in $n$ variables $x_1-p_1, \cdots, x_n-p_n$, where $p=(p_1,\cdots,p_n) \in \r^{n}$.

\section{Representation of nonnegative Morse polynomial functions}
 In \cite[Theorem 1.6.4]{Mar08} the author showed that if a real polynomial  in one variable (resp.  two variables) which is  nonnegative in a neighborhood of $0\in \r$ (resp.  $(0,0)\in \r^{2}$),  then $f\in \r[[x_1]]^{2}$ (resp.  $f\in \sum \r[[x_1,x_2]]^{2}$). 
Moreover, for $n\geq 3$,  there exists always a polynomial which is nonnegative on $\r^{n}$ but  does not belong to $\sum \r[[X]]^{2}$.  However, for a Morse polynomial function, we have a nice representation.

\lm \label{lm-basic} Let $f: \r^{n} \mtn \r $ be a Morse polynomial function. Assume that $f(x)\geq 0 $ for every $x$ in a neighborhood $U$ of $0\in \r^{n}$, and $f^{-1}(0)=\{0\}$. Then $f\in \sum\r[[X]]^{2}$.
\elm 
\pf  It follows from the assumption that $0\in \r^{n}$ is an isolated minimal point of $f$.   The minimality of the local minimum $0$ of $f$  implies that the Hessian matrix  $D^{2}f(0)$  of $f$ at $0$ is positive semidefinite, i.e. all of its eigenvalues are nonnegative.  
On the other hand, since   $0\in \r^{n}$ is a critical point of $f$ who is Morse,    $0$ is non-degenerate, i.e. the Hessian matrix $D^{2}f(0)$   is invertible. Therefore the matrix  $D^{2}f(0)$ has no zero eigenvalues, i.e. all eigenvalues of $D^{2}f(0)$  are positive. Thus the Hessian matrix $D^{2}f(0)$ of $f$ at $0$ is positive definite.   \\
Then by a linear  change of coordinates in a neighborhood of $0\in \r^{n}$ we may assume that in a neighborhood of $0\in \r^{n}$ the polynomial $f$ is expressed in  the following form: 
$$ f= x_1^{2} + \cdots + x_n^{2} + g,$$
where the order of $g$ is greater than or equal to $3$. For a monomial $a_\alpha x_1^{\alpha_1}\cdots x_n^{\alpha_n}$ in $g$ such that $\sum \alpha_i \geq 3$ and  there exists $i \in \{1,\cdots, n\}$ such that $\alpha_i \geq 2$, we have 
\begin{equation}\label{pt-1}
 \epsilon_i x_i^{2} + a_\alpha x_1^{\alpha_1}\cdots x_n^{\alpha_n} = x_i^{2}(\epsilon_i + a_\alpha x_1^{\alpha_1}\cdots x_i^{\alpha_i-2}\cdots x_n^{\alpha_n}) \in \r[[X]]^{2}, 
\end{equation}
where $0<\epsilon_i \ll 1$. Note that the inclusion in (\ref{pt-1}) follows from the fact that for  $g\in \r[[X]]$ with $g(0)>0$ we have $g\in \r[[X]]^{2}$ (cf. \cite[Proposition 1.6.2]{Mar08}).  Therefore,  by renumbering the indices if necessary, it suffices to prove that 
$$ h(x_1,\cdots,x_m):=x_1^{2} + \cdots + x_m^{2} + a x_1 \cdots x_m \in \sum \r[[X]]^{2}, $$
where $a\in \r$ and $m\geq 3$. In fact, for each $i=1,\cdots, m$, denote $u_i:=\prod_{j\not =i} x_j$. Then 
$$ h=\sum_{i=1}^{m} \Big(\ps{1}{2}x_i + \ps{a}{m}u_i\Big)^{2}  + \ps{3}{4} \sum_{i=1}^{m} x_i^{2}  - \ps{a^{2}}{m^{2}} \sum_{i=1}^{m} u_i^{2}.$$
Note that for any $b\in \r$ and for any $i\not =j$, similar to the argument  shown above, we have 
$$ x_i^{2} + b u_j^{2} = x_i^{2}(1+b \ps{u_j^{2}}{x_i^{2}}) \in \r[[X]]^{2}.$$
Then $h \in \sum \r[[X]]^{2}$. The proof is complete.
\epf 

\thr[\textbf{Scheiderer's Positivstellensatz for Morse polynomials}]  \label{thr-PosMorse}Let $G=\{g_1,\cdots, g_m\}\subseteq \r[X]$, and $K_G$ be the basic closed semi-algebraic set generated by $G$.  Let $f: \r^{n} \mtn \r$ be a Morse polynomial function. Assume the following conditions hold true:
\ite 
\item[(1)] $K_G$ is compact (resp. $M_G$ is Archimedean);
\item[(2)]  $f\geq 0$ on $K_G$, and  $f$ has only finitely many zeros $p_1,\cdots, p_r$ in $K_G$, each lying in the interior of $K_G$.
\hite
 Then $f \in T_G$ (resp. $f\in M_G$).
\ethr
\pf   Since each $p_i$ is an isolated minimal point of $f$, it follows from Lemma \ref{lm-basic} that $f \in \sum \r[[X-p_i]]^{2}$.  On the other hand,  for every $i=1,\cdots,r, $   $p_i$ is an interior point of $K_G$,  hence the condition (4') in Theorem \ref{local-global2} is fulfilled and by  Lemma \ref{lm-interior} below,   we have 
$$ \widehat{T}_{p_i} = \widehat{M}_{p_i} = \sum \r[[X-p_i]]^{2}. $$
The theorem  now follows from Theorem \ref{local-global1} and Theorem \ref{local-global2}.
\epf 
\lm[{\cite[Lemma 2.1]{Le}}] \label{lm-interior} If $p\in K_G$ is an interior point then 
 $$\sum \r[[X-p]]^2=\widehat{T}_{p}=\widehat{M}_{p}.$$
\elm 
\pf 
 Since $p$ is an interior point of $K_G$, we have $g_i(p)>0$, for all $i=1,\cdots, m$. Then   $g_i\in \r[[X-p]]^2$ for all $i=1,\cdots,m$ (cf. \cite[Proposition 1.6.2]{Mar08}). It follows that  $\widehat{M}_p \subseteq  \widehat{T}_p\subseteq  \sum \r[[X-p]]^2$.  It is clear that $\sum \r[[X-p]]^2 \subseteq \widehat{M}_p$. Thus we have the equalities.
\epf 

\rem \label{rem-interior} \rm (1) The assumption that each zero of $f$  belongs to the interior of $K_G$ in Theorem \ref{thr-PosMorse} is necessary. Indeed,  let us consider again $G=\{(1-x^2)^3\}\subseteq \r[x]$ and $f=1-x^2 \in \r[x]$.  We see that  $K_G=[-1,1]$ is compact,   $f$ is a Morse function (it has a critical point at $0\in \r$ and the second derivative of $f$ at $0$ is equal to $-2$ which is non-zero), $f\geq 0$ on $K_G$, and $f$ has only two zeros in $K_G$. However, $f \not \in T_G=M_G$. In this case, note that the zeros of $f$ belong to the boundary of $K_G$. \\
(2) The space of Morse functions on  $\r^n$ is a \textit{dense} subset of the space of all smooth  functions on $\r^n$ in the uniform topology (cf. \cite[Theorem 5.27]{BH04}). Therefore Theorem \ref{thr-PosMorse} holds for  a generic class of polynomial functions on $\r^n$. \\
(3) In \cite[Theorem 2.33]{La10}, J.-B. Lasserre has given a similar Positivstellensatz for the case where  $f$ is \textit{stricly convex } and  $g_j$ is  \textit{concave} for every  $j=1,\cdots, m$.   

\erem

\section{Representation of Morse polynomial functions on non-degenerate algebraic sets}
As we have seen in the previous section, the compactness of the semi-algebraic sets $K_G$ generated by a finite subset $G$ in $\r[X]$ is very important for the representation of a nonnegative polynomial function on $K_G$. In this section we give a good class of polynomials  in $\r[X]$ for which the algebraic sets defined by them are compact. For this purpose, we introduce the   \textit{non-degeneracy conditions} which was studied by C. Bivi\`{a}-Ausina \cite{Bi07}.

\df[{\cite{Bi07}, \cite{DHP14}}] \rm  A subset $\widetilde{\Gamma} \subseteq \r_{+}^{n}$ is said to be a \textit{Newton polyhedron at infinity}, or a \emph{global Newton polyhedron},  if there exists some finite  subset $A$ of $\mathbb Z_{+}^{n}$ such that $\widetilde{\Gamma}$ is  equal to the convex hull of $A\cup \{0\}$ in $\r^{n}$.  $\widetilde{\Gamma}$ is said to be \textit{convenient}, if it  intersects each coordinate axis in a point different from the origin. 

For $w\in \r^n$, denote 
$$ m(w,\widetilde{\Gamma}):=\max\{\seq{w,\alpha}| \alpha \in \widetilde{\Gamma}\}; $$
$$ \Delta(w,\widetilde{\Gamma}):=\{\alpha\in \widetilde{\Gamma}| \seq{w,\alpha} = m(w,\widetilde{\Gamma})\}. $$
A set $\Delta \subseteq \r^n$ is called a face of $\widetilde{\Gamma}$ if there exists some $w\in \r^n$ such that $\Delta=\Delta(w,\widetilde{\Gamma})$. In this case, the face $\Delta(w, \widetilde{\Gamma})$ is said to be \textit{supported} by $w$.  

 Let $f=\sum_\alpha f_\alpha X^\alpha\in \r[X]$ be a polynomial and $w\in \r^n$. The set 
$ supp(f):=\{\alpha \in \n^n| f_\alpha \not =0\} $
is called the \textit{support} of $f$.  Denote 
$$ m(w,f):=\max \{\seq{w,k}| k\in supp(f)\}; $$
$$ \Delta(w,f):=\{k\in supp(f)| \seq{w,k} = m(w,f) \}.$$
The convex hull in $\r_{+}^n$ of the set $supp(f)\cup \{0\}$ is called the \textit{Newton polyhedron at infinity} of $f$ and denoted by $\widetilde{\Gamma}(f)$. We say that $f$ is \textit{convenient}  if $\widetilde{\Gamma}(f)$  is convenient. \\
The polynomial $f_w:=f_{\Delta(w,f)}:=\sum_{\alpha \in \Delta(w,f)} f_\alpha X^\alpha$ is called the \textit{principal part of $f$ at infinity} with respect to $w$ (or $\Delta(w,f)$).  For a finite subset $W$ of $\r^n$, the \textit{principal part of $f$ wih respect to $W$ at infinity} is defined to be the polynomial $f_W:=\sum_{\alpha \in \cap_{w\in W} \Delta(w,f)} f_\alpha X^\alpha$. If $\cap_{w\in W} \Delta(w,f) = \emptyset$, we set $f_W=0$.

Let $F=(f_1,\cdots,f_m): \r^n \mtn \r^m$ be a polynomial map. Then the convex hull of $\widetilde{\Gamma}(f_1) \cup \cdots \widetilde{\Gamma}(f_m)$ is called the \textit{Newton polyhedron at infinity of $F$} and denoted by $\widetilde{\Gamma}(F)$. For $w\in \r^n$, the \textit{principal part of $F$ with respect to $w$ at infinity }  is defined to be  the polynomial map
$$ F_w:=((f_1)_w,\cdots, (f_m)_w). $$
\edf 

\df[{\cite{Bi07}}] \rm Let $F=(f_1,\cdots,f_m): \r^n \mtn \r^m$ be a polynomial map. Denote 
$$ \r_0^n:=\{w=(w_1,\cdots,w_n)\in \r^n| \max_{i=1,\dots,n} w_i >0\}. $$
We say that $F$ is \textit{non-degenerate at infinity} if and only if for any $w\in \r_0^n$, the system of equations 
$$ (f_1)_w(x) = \cdots = (f_m)_w(x) =0 $$
has no solutions in $(\r\tru \{0\})^n$.
\edf 
\rem[{\cite{Bi07}}] \rm  \rm (1) If some component of $F$ is a monomial, then $F$ is automatically non-degenerate at infinity.\\
(2) Let $F=(f_1,\cdots,f_m): \r^2 \mtn \r^m$ such that $\widetilde{\Gamma}(f_i)$ is convenient for every $i=1,\cdots,m$. Let $\widetilde{\Gamma}_\vc(f_i)$ denote the \textit{Newton boundary at infinity} of $f_i$, i.e. the union of all faces of $\widetilde{\Gamma}(f_i)$ which do not passing through the origin.  Then  $F$ is non-degenerate at infinity if either some component $f_i$ is a monomial or the polygons of the family $\{\widetilde{\Gamma}_\vc(f_1), \cdots, \widetilde{\Gamma}_\vc(f_m)\}$ verify that no segment of  $\widetilde{\Gamma}_\vc(f_i)$ is parallel to some segment of $\widetilde{\Gamma}_\vc(f_j)$ for all $i,j\in \{1,\cdots,m\}$, $i\not = j$. 
\erem 
\thr[{\cite[Theorem 3.8]{Bi07}}] \label{thr-compact1} Let $F=(f_1,\cdots,f_m): \r^n \mtn \r^m $ be a polynomial map such that  
$f_i$ is convenient for all $i=1,\cdots,m$. If $F$ is non-degenerate at infinity, then $F^{-1}(0)$ is compact. 
\ethr 

The algebraic set $K_G:=\{x\in \r^{n}| g_1(x)=\cdots =g_m(x)=0\}$, $g_i\in \r[X] $ for all $i=1,\cdots,m$,  is called \emph{non-degenerate at infinity} if the polynomial map $(g_1,\cdots,g_m): \r^{n} \mtn \r^{m}$ is non-degenerate at infinity. Then we have the following special case of Theorem \ref{thr-PosMorse}.

\coro \label{coro-Morsecompact1} Let $G=\{g_1,\cdots, g_m\}$ be a finite subset of $\r[X]$ and $K_G:=\{x\in \r^{n}| g_i(x)=0, i=1,\cdots,m\}$ the algebraic set defined by $G$.  Let $f: \r^{n} \mtn \r$ be a Morse polynomial function.   Assume the following conditions hold true:
\ite 
\item[(1)]   $K_G$ is non-degenerate at infinity and each $g_i$ is convenient;
\item[(2)]  $f\geq 0$ on $K_G$, and  $f$ has only finitely many zeros $p_1,\cdots, p_r$ in $K_G$, each lying in the interior of $K_G$.
\hite
 Then $f \in T_G$. 
\ecoro 
\pf The proof follows from Theorem \ref{thr-compact1} and Theorem \ref{thr-PosMorse}.
\epf 

\rem \rm In Theorem \ref{thr-compact1} we need the convenience of each component $f_i$ of the polynomial map $F=(f_1,\cdots,f_m)$ for $F^{-1}(0)$ to be compact. In the following we introduce another condition of non-degeneracy for which the assumption on the convenience of each $f_i$ can be relaxed. 
\erem 

\df[\cite{Bi07}] \rm  Let $f:(\r^{n},0) \mtn (\r,0)$ be a real analytic function.  Suppose that the Taylor expansion of $f$ around the origin
is given by the expression $f=\sum_{\alpha} f_\alpha X^{\alpha}$. The set $supp(f):=\{\alpha \in \mathbb N^{n}| f_\alpha \not =0\}$ is called the \textit{support }of $f$.  For a vector $v\in \r^n_{+}$, denote 
$$ l(v,f):=\min \{\seq{v,\alpha}| \alpha \in supp(f)\}. $$
For a finite set $V$ of $\r^n_{+}$, the \textit{local principal part of $f$ with respect to $V$} is defined to be the polynomial 
$$ f_V:=\sum_{\seq{\alpha, v} = l(v,f), \vm v\in V } f_\alpha X^\alpha. $$
If  no such terms exist we define $f_V=0$. 

The \textit{local Newton polyhedron of $f$}, denoted by $\Gamma(f)$, is the convex hull of the set 
$$ \bigcup_{\alpha \in supp(f)}\{\alpha+ \r_{+}^{n}\}. $$
A subset $\Gamma $ of $\r_{+}^{n}$ is said to be a \textit{local Newton polyhedron} if there exists some real analytic function $f$ such that $\Gamma = \Gamma(f)$. 

Let $\Gamma$ be a local Newton polyhedron in $\r_+^n$. For $v\in \r_+^n$, we define 
$$ l(v,\Gamma):=\min\{\seq{v,\alpha}| \alpha \in \Gamma\}; $$
$$ \Delta(v,\Gamma):=\{\alpha\in \Gamma| \seq{v,\alpha} = l(w,\Gamma)\}. $$
A set $\Delta\subseteq \r_{+}^n$ is called a \textit{face} of $\Gamma$ if there exists some $v\in \r_{+}^n$ such that $\Delta=\Delta(v,\Gamma)$. Then we say that the vector $v$  \textit{supports}  the face $\Delta$.  \\
A vector $w\in \z^n$ is called \textit{primitive}  if $w\not = 0$ and it has smallest length among all vectors in $\z^n$ of the form $\lambda w, \lambda >0$.  Denote by $\mathcal{F}(\Gamma)$ the family of primitive vectors supporting some face of $\Gamma$ of dimension $n-1$.  
\edf 

\df[\cite{Bi07}] \rm Let $\Gamma$ be a local Newton polyhedron in $\r_{+}^n$.   Let $f=(f_1,\cdots,f_m): (\r^n,0) \mtn (\r^m,0)$ be an analytic map germ. $f$ is said to be \textit{adapted to $\Gamma$} if for all $V\subseteq \mathcal{F}(\Gamma)$ such that $\cap_{v\in V}\Delta(v,\Gamma)$ is a compact face of $\Gamma$, the system of equations
$$ (f_1)_V (x) = \cdots = (f_m)_V(x) = 0 $$
has no solutions in $(\r\tru \{0\})^n$.
\edf 

\df[\cite{Bi07}]  \rm For $I\subseteq \{1,\cdots,n\}$, denote 
$$ \r_I^n = \{x=(x_1,\cdots,x_n)\in \r^n| x_i=0 \mbox{ for all } i\in I\}. $$
If $I=\emptyset $ then it is clear that  $\r_I^n = \r^n$. \\
For a polynomial $f=\sum_\alpha f_\alpha X^\alpha \in \r[X]$, denote 
$$ f_I:=\sum_{\alpha \in \r_I^n} f_\alpha X^\alpha. $$
If $supp(f) \cap \r_I^n=\emptyset$, we set $f_I=0$.  We regard  $f_I$ as a polynomial in  on the variables $x_i$ such that $i\not \in I$, i.e. $f_I$ can be regarded as  the function  $f_I: \r^{n-|I|} \mtn \r$.  For a polynomial map $F=(f_1,\cdots,f_m): \r^n \mtn \r$, $F_I$ denotes the map $((f_1)_I,\cdots,(f_m)_I): \r^{n-|I|} \mtn \r$.   

Let $\widetilde{\Gamma}$ be a \textit{fixed convenient Newton polyhedron at infinity} in $\r^{n}$ and $I\subseteq \{1,\cdots,n\}$. Denote  by $(\widetilde{\Gamma})_I$ the image of the intersection $\widetilde{\Gamma}\cap \r_{I}^n$ in $\r^{n-|I|}$. Set 
$$ M:=\max_{\alpha=(\alpha_1,\cdots,\alpha_n) \in \widetilde{\Gamma}} \{ |\alpha|:=\alpha_1+\cdots + \alpha_n\}.$$
Let $V_{\widetilde{\Gamma}}$ denote the set of all vertices of $\widetilde{\Gamma}$, and $\rho:=\sum_{\alpha \in V_{\widetilde{\Gamma}}} X^{\alpha}$. For any polynomial $h=\sum_\alpha h_\alpha X^{\alpha}\in \r[X]$, denote 
$$ G_M(h):=\sum_{\alpha} h_\alpha X^{\alpha} \|x\|^{2(M-|\alpha|)}.$$
Then we define the convenient local Newton polyhedron associated to the global Newton polyhedron $\widetilde{\Gamma}$:
$$ \bold{G}(\widetilde{\Gamma}) :=\Gamma(G_M(\rho)).$$
\edf
\df[\cite{Bi07}]  \rm  Let $F=(f_1,\cdots,f_m): \r^n \mtn \r^m$ be a polynomial map. We say that $F$ is \textit{globally adapted  to $\widetilde{\Gamma}$} (or, \textit{g-adapted to $\widetilde{\Gamma}$}) if for any $W\subseteq \{\bold{w}(v)| v\in \mathcal{F}(\bold{G}(\widetilde{\Gamma}))\}$ such that $\cap_{w\in W}\Delta(w,\widetilde{\Gamma})$ is a face of $\widetilde{\Gamma}$ not containing the origin,  the system of equations
$$ (f_1)_W (x) = \cdots = (f_m)_W(x) = 0 $$
has no solutions in $(\r\tru \{0\})^n$. 
Here, for a vector $v=(v_1,\cdots,v_n)\in \r^n$, $\bold{w}(v):=2\bold{c} \min_{i} v_i -v$, where $\bold{c}:=\bold{e}_1 +\cdots + \bold{e}_n = (1,\cdots,1) \in \r^n$. 

We say that $F$ is \textit{strongly g-adapted} to $\widetilde{\Gamma}$ if for any $I\subseteq \{1,\cdots,n\}$, $|I| \not = n$, the map $F_I: \r^{n-|I|} \mtn \r^m$ is g-adapted to the Newton polyhedron $(\widetilde{\Gamma})_I$.

It follows from the above definition that if $F$ is strongly g-adapted to a given convenient Newton polyhedron at infinity then $\widetilde{\Gamma}(F)$ is convenient.
\edf 

\thr[{\cite[Theorem 5.9]{Bi07}}] \label{thr-compact2}  Let $\widetilde{\Gamma}$ be a  convenient   Newton polyhedron at infinity.     Let $F=(f_1,\cdots,f_m): \r^n \mtn \r^m$ be a polynomial map with degree $d:=\max\{\deg(f_1),\cdots,\deg(f_m)\}$ such that $M\geq d$.  If $F$ is strongly g-adapted to $\widetilde{\Gamma}$ then $F^{-1}(0)$ is compact.
\ethr 

\coro \label{coro-Morsecompact2} Let $G=\{g_1,\cdots, g_m\}$ be a finite subset of $\r[X]$ and $K_G:=\{x\in \r^{n}| g_i(x)=0, i=1,\cdots,m\}$ the algebraic set defined by $G$.   Let $\widetilde{\Gamma}$ be a  convenient   Newton polyhedron at infinity.   Let $f: \r^{n} \mtn \r$ be a Morse polynomial function.   Assume the following conditions hold true:
\ite 
\item[(1)]   The polynomial map $(g_1,\cdots,g_m): \r^n \mtn \r^m$ has degree $\leq M$ and  strongly g-adapted to $\widetilde{\Gamma}$; 
\item[(2)]  $f\geq 0$ on $K_G$, and  $f$ has only finitely many zeros $p_1,\cdots, p_r$ in $K_G$, each lying in the interior of $K_G$.
\hite
 Then $f \in T_G$. 
\ecoro 
\pf  The proof follows from Theorem \ref{thr-compact2} and Theorem \ref{thr-PosMorse}.
\epf

\section{Applications in polynomial Optimization}
\subsection{Unconstrained polynomial Optimization}
In this section we consider the global  optimization problem 
\begin{equation}\label{pt-global}
f^*=\min_{x\in \r^n} f(x),
\end{equation}
where $f\in \r[X]$ be a polynomial in $n$ variables $x_1,\cdots, x_n$. 

It is well-known (cf. \cite{NDS06}) that if the gradient ideal $I_{grad}(f)$ is radical and if $f$ attains its minimum value $f^*$ on $\r^n$, then $f-f^*$ is SOS modulo $I_{grad}(f)$.  In general we have $f-f^*$ is SOS modulo the radical $\can{I_{grad}(f)}$ of the gradient ideal $I_{grad}(f)$ (cf. \cite{NDS06}).  However, for Morse polynomial functions we have a nice representation of $f-f^*$. 

\pro \label{pro-sosgrad1} Let $f: \r^n \mtn \r$ be a Morse polynomial function. Assume that $f$ achieves a  minimum value $f^*$ on $\r^n$. Then 
$$f-f^* \in \sum \r[X]^2 + I_{grad}(f),$$
where $I_{grad}(f) = \seq{\ps{\partial f}{\partial x_1}, \cdots, \ps{\partial f}{\partial x_n}}$ denotes the gradient ideal of $f$. 
\epro 
\pf  Let $x^{*}\in \r^{n}$ be a global minimizer of $f$ on $\r^{n}$. Then $x^{*}$ is a critical point of $f$, therefore  the Hessian matrix $D^{2}f(x^{*})$ is invertible because $f$ is Morse. Moreover, $D^{2}f(x^{*})$ is positive semidefinite because $x^{*}$ is a global minimizer. It follows that $D^{2}f(x^{*})$ is positive definite. Now apply  \cite[Theorem 2.1]{Mar09}, we have  
$$f-f^* \in \sum \r[X]^2 + I_{grad}(f).$$
\epf 
The following result gives  degree bounds to accompany Proposition \ref{pro-sosgrad1}. 
\pro \label{pro-degreebound} Given a positive integer $d$. Then there exists a positive integer $l$ such that for each Morse polynomial function  $f: \r^n \mtn \r$  of degree $\leq d$, if  $f$ achieves a  minimum value $f^*$ on $\r^n$, then 
$$f-f^*=\sigma + \sum_{i=1}^{n}h_i \ps{\partial f}{\partial x_i},$$
where $\sigma \in \sum \r[X]^{2}$ and $h_1,\cdots, h_n \in \r[X]$, have degree bounded by $l$.
\epro 
\pf  Similar to the proof of Proposition \ref{pro-sosgrad1}, if $f$ is Morse and $x^{*}$ is a global minimizer of $f$ on $\r^{n}$ then the Hessian matrix $D^{2}f(x^{*})$ is positive definite. Then the proposition follows from  \cite[Corollary 2.4]{Mar09}.
\epf 

Let $\r[X]_m$ denote the $\binom{n+m}{m}$-dimensional vector space of polynomials of degree at most $m$. Since the gradient is zero at   global minimizers, we consider the SOS relaxation:
\begin{align}\label{pt-SDPgrad}
f^*_{N,grad}&:=\max   \gamma \\
&  \mbox{ subject to } f-\gamma-\sum_{i=1}^n \phi_i\ps{\partial f}{\partial x_i} \in \sum \r[X]^2  \mbox{ and } \phi_i \in \r[X]_{2N-d+1}.\nonumber
\end{align}
Here $d$ is the degree of the polynomial $f\in \r[X]$, and $N$ is an integer to be chosen by the user. 

It is well-known (cf. \cite{La01}, \cite{Lau09}, \cite{Mar03}, \cite{NDS06}, \cite{PS03}) that the problem (\ref{pt-SDPgrad}) can be translated into an SDP.   Moreover,  $f^{*}_{N,grad}$ is a lower bound for $f^{*}$, and  the lower bound gets better as $N$ increases:
$$ \cdots \leq f^{*}_{N-1,grad} \leq  f^{*}_{N,grad} \leq f^{*}_{N+1,grad}\leq \cdots \leq f^{*}.$$
In the following we apply Proposition  \ref{pro-sosgrad1} to show the finite convergence of the relaxation given above in the case where $f$ is a Morse polynomial function.

\thr \label{thr-sosconvergence1} Let $f: \r^n \mtn \r$ be a Morse polynomial function. Assume that $f$ achieves a  minimum value $f^*$ on $\r^n$. Then  there exists an integer $N$ such that $f^{*}_{N,grad}=f^{*}$.
\ethr
\pf  It follows from Proposition \ref{pro-sosgrad1} that $f-f^{*}  $ is SOS modulo $I_{grad}(f)$. Then   by Proposition \ref{pro-degreebound}, there  exists some positive integer $N$ such that $  f^{*}_{N,grad}\geq f^{*}$. Moreover, we have always that $ f^{*}_{N,grad}\leq f^{*}$. Hence $f^{*}_{N,grad}=f^{*}$. 
\epf 

\rem \rm (1) The assumption that $f$ achieves  a  minimum value $f^*$ on $\r^n$ is necessary. Indeed, let us consider the polynomial $f(x) =x^{3}$ in one variable. It is clear that $f^{*} = -\vc$ on $\r$.   Moreover, we have 
$$ f(x) = \ps{x}{3} f'(x), $$
hence $f$ belongs to its gradient ideal $I_{grad}(f)=\seq{f'}$. Therefore for every $N\geq 1$ we have $f^{*}_{N,grad} = 0>f^{*}$. \\
(2) There is a generic class of polynomials which achieve their minimum values on $\r^{n}$. For example, in \cite[Theorem 1.1]{DHP14} the authors showed that if $f\in \r[X]$ is bounded from below,  convenient\footnotemark \footnotetext{The polynomial function $f: \r^n \mtn \r$ is said to be \textit{convenient} if its \textit{Newton polyhedron at infinity} $\widetilde{\Gamma}(f)$  intersects each coordinate axis in a point different from the origin, that is, if for any $i\in \{1,\cdots,n\}$ there exists some integer $m_i>0$  such that $m_i\bold{e}_i \in \widetilde{\Gamma}(f)$. Here $\{\bold{e}_1,\cdots,\bold{e}_n\}$ denotes the canonical basis in $\r^n$.  } and (Khovanskii) non-degenerate at infinity\footnotemark, \footnotetext{$f$ is called \emph{(Khovanskii) non-degenerate at infinity} if for any face $\Delta$ of $\widetilde{\Gamma}(f)$ which does not contain the origin $0\in \r^n$, the system of equations 
$$ f_\Delta=x_1\ps{\partial f_\Delta}{\partial x_1}=\cdots =  x_n\ps{\partial f_\Delta}{\partial x_n}=0$$
has no solution in $(\r\tru \{0\})^n$. Here for $f=\sum_{\alpha} f_\alpha X^\alpha \in \r[X]$,  $f_\Delta:=\sum_{\alpha \in \Delta} f_\alpha X^\alpha$ denotes the \textit{principal part at infinity } of $f$ with respect to $\Delta$. } then $f$ attains its minimum value $f^{*}$ on $\r^{n}$.  Therefore we have the following consequence of this fact and Theorem
 \ref{thr-sosconvergence1}.
\erem 
\coro \label{coro-sosnondegenerate} Let $f: \r^n \mtn \r$ be a Morse polynomial function which is bounded from below, convenient and (Khovanskii) non-degenerate at infinity.  Then $f$ achieves its minimum value $f^*$ on $\r^n$, moreover,  there exists an integer $N$ such that $f^{*}_{N,grad}=f^{*}$. 
\ecoro 

\subsection{Constrained polynomial optimization} Let $G=\{g_1,\cdots, g_m\}$ be a finite subset of $\r[X]$ and $K_G=\{x\in \r^n| g_i(x) \geq 0, \forall i=1,\cdots, m\}$ the basic closed semi-algebraic set generated by $G$. In this section we consider the following optimization problem 
\begin{equation}\label{pt-constrained}
f^*:=\min_{x\in K_G} f(x),
\end{equation}
where $f\in \r[X]$ be a polynomial in $n $ variables $x_1,\cdots, x_n$. The \textit{KKT system} associated to this optimization problem is 
\begin{align} \label{kkt-system}
\nabla f - \sum_{j=1}^m \lambda_j \nabla g_j = 0\\
g_j \geq 0, \quad \lambda_j g_j \geq 0, \quad j=1,\cdots, m\nonumber
\end{align}
where the variables  $\lambda:=(\lambda_1,\cdots,\lambda_m)$ are called \textit{Lagrange multipliers} and $\nabla  f$ denotes the vector of partial derivatives of $f$. A point is called a \textit{KKT point } if  the KKT system holds at this point. Under certain regularity conditions, for example if the gradients $\nabla g_i $ of the $g_i$'s  are linearly independent (cf. \cite{NW}),  each  global minimizer of $f$ on $K_G$ is a KKT point. 

For each $i=1,\cdots,n$, denote  $ L_i:=\ps{\partial f}{\partial x_i} - \sum_{j=1}^m \lambda_j \ps{\partial g_j}{\partial x_i},  i=1,\cdots, n. $ We define the \textit{KKT ideal } $I_{KKT}$,  the \textit{KKT varieties }, the  \emph{KKT preordering}  and the \emph{KKT quadratic module} associated to the KKT system (\ref{kkt-system}) as follows.
$$ I_{KKT}:=\seq{L_1,\cdots,L_n, \lambda_1g_1,\cdots,\lambda_mg_m}; $$
$$ V_{KKT}:=\{(x,\lambda) \in \c^n\times \c^m| g(x)=0 \mbox{ for all } g\in I_{KKT}\}; $$
$$ V^\r_{KKT}:= \{(x,\lambda) \in \r^n\times \r^m| g(x)=0 \mbox{ for all } g\in I_{KKT}\};$$
$$ T_{KKT}:=T_G+I_{KKT}; $$
$$ M_{KKT}:= M_G + I_{KKT}. $$
Let $f^*_{KKT}$ be the global minimum of $f$ over the KKT system defined by (\ref{kkt-system}). Assume the KKT system
holds at at least one global minimizer. Then $f^*=f^*_{KKT}$ (cf. \cite{NDP07}). Therefore we have 
\begin{equation}\label{min-kkt}
f^*=\min_{x\in V^\r_{KKT}\cap K_G} f(x)
\end{equation}
provided that the KKT system  holds at at least one global minimizer. 

It is well-known (cf. \cite{NDP07}),  that  if $I_{KKT}$ is zero-dimensional (i.e. $V_{KKT}$ is a finite set) and radical, then $f-f^* \in M_{KKT}$.  Moreover, if $I_{KKT}$ is radical, then $f-f^* \in T_{KKT}$ (cf. \cite{NDP07}).  For Morse polynomial functions we have 
\pro \label{pro-kkt-representation} Let $G=\{g_1,\cdots, g_m\}$ be a finite subset of $\r[X]$ such that $K_G$ is compact (resp. $M_G$ is Archimedean).  Let $f:\r^n \mtn \r$ be a Morse polynomial function.   Assume that $V^\r_{KKT}\cap K_G$ is finite and contained in the interior of $K_G$.  Then $f-f^* \in T_{G}$ (resp. $f-f^* \in M_G$).
\epro 
\pf  It is obvious that $f-f^*$ is a Morse polynomial function. Moreover,  $f-f^* \geq 0$ on $K_G$, and by assumption, $f-f^*$ vanishes at only finitely many points in the interior of $K_G$. Then it follows from Theorem \ref{thr-PosMorse}  that $f-f^* \in T_{G}$ (resp. $f-f^* \in M_G$).
\epf 

To give more applications in polynomial optimization we need to recall some notations (cf. \cite{La01}). Let 
$$
\bold{v}_m(\bold{x}):=\big(1,x_1,\cdots,x_n, x_1^2,x_1x_2,\cdots,x_1x_n,x_2x_3,\cdots,x_n^2,\cdots, x_1^m,\cdots,x_n^m\big)
$$
denotes the canonical basis of the real vector space $\r[X]_{m}$ of real polynomials of degree at most $m$, and let $s(m):=\binom{n+m}{m}$ be the dimension of this  vector space.   If $f$ is a polynomial of degree at most $m$ we may write 
$$ f=\sum_\alpha f_\alpha X^\alpha = \seq{\bold{f}, \bold{v}_m(\bold{x})},  \mbox{ where } X^\alpha:=x_1^{\alpha_1}\cdots x_n^{\alpha_n}, \sum_{i=1}^n \alpha_i \leq m, $$
and $\bold{f}:=\{f_\alpha\}\in \r^{s(m)}$ denotes the vector of coefficients of $f$ in the basis $\bold{v}_m(\bold{x})$. 

Given an $s(2m)$-vector $\bold{y}:=\{y_\alpha \}$ with first element $y_{0,\cdots,0} = 1$, let $M_m(\bold{y})$ be
the \textit{moment matrix} of dimension $s(m)$, with rows and columns labeled by the basis $\bold{v}_m(\bold{x})$. 

Let  $f\in \r[X]_{m}$ with coefficient vector $\bold{f}\in \r^{s(m)}$. If the entry $(i,j)$ of the matrix $M_m(\bold{y})$ is $y_\beta$, let $\beta(i,j)$ denote the subscript $\beta$ of $y_\beta$. We define the matrix $M_m(f\bold{y})$ by
$$ M_m(f\bold{y})(i,j):=\sum_\alpha f_\alpha y_{\beta(i,j)+\alpha}. $$

Now let $G=\{g_1,\cdots,g_m\}$ be a finite subset of $\r[X]$, with each  $g_i$ is a  polynomial of degree at most $w_i$.  Let  $f\in \r[X]_{m}$ with coefficient vector $\bold{f}=\{f_\alpha\}\in \r^{s(m)}$. For every $i=1,\cdots,m$, let $\tilde{w_i}=\lceil w_i/2 \rceil$ be  the smallest integer larger than $w_i/2$, and with $N\geq \lceil m/2 \rceil$ and $N\geq \max_i \tilde{w_i}$, consider the convex LMI problem 
$$ \mathbb Q^N_{G}
\begin{cases} 
\inf_{\bold{y}} \sum_{\alpha} f_\alpha y_\alpha, & \\
  \quad M_N(\bold{y}) &\succcurlyeq 0, \\
M_{N-\tilde{w_i}} (g_i\bold{y}) &\succcurlyeq 0, \quad i=1,\cdots,m.  
\end{cases}
$$

\thr[{\cite[Theorem 4.2]{La01}}] \label{thr-Lasserre} Let $G=\{g_1,\cdots,g_m\}$ be a finite subset of $\r[X]$ and $K_G$ the basic closed semi-algebraic generated by $G$. Assume $M_G$ is Archimedean. Let $f\in \r[X]$ be a polynomial of degree $m$. If there exist a  polynomial $q\in \sum \r[X]^2$ of degree at most $2N$ and polynomials  $ t_i \in \sum \r[X]^2$ of degree at most $2N-w_i$, $i=1,\cdots, m$, such that 
$$ f-f^*=q + \sum_{i=1}^m t_ig_i, $$
then $\min \mathbb Q^N_G = f^*$, and the vector 
$$ \bold{y}^*:=\big(x_1^*,\cdots,x_n^*,(x_1^*)^2,\cdots, x_1^*x_2^*,\cdots, (x_1^*)^{2N},\cdots, (x_n^*)^{2N}\big) $$ 
is a global minimizer of $\mathbb Q^N_G$.
\ethr 
Combining  Proposition \ref{pro-kkt-representation} and Theorem  \ref{thr-Lasserre}, we have the following result.

\coro \label{coro-momentrepresentation} Let $G=\{g_1,\cdots, g_m\}$ be a finite subset of $\r[X]$ such that   $M_G$ is Archimedean.  Let $f:\r^n \mtn \r$ be a Morse polynomial function.    Assume that $V^\r_{KKT}\cap K_G$ is finite and contained in the interior of $K_G$.  Then there exists a positive integer $N$ such that $\min \mathbb Q^N_G = f^*$. Moreover,  if $x^{*} \in  V^\r_{KKT}\cap K_G$  is a global minimizer of $f$ on $K_G$, then the vector 
$$ \bold{y}^*:=\big(x_1^*,\cdots,x_n^*,(x_1^*)^2,\cdots, x_1^*x_2^*,\cdots, (x_1^*)^{2N},\cdots, (x_n^*)^{2N}\big) $$ 
is a global minimizer of $\mathbb Q^N_G$.
\ecoro 

\hspace{-0.6cm} \textbf{Acknowledgements} ~ The author would like to express his gratitude to Professor H\`{a} Huy Vui  for his  valuable discussions on Morse theory and polynomial optimization.   

 The original verson of this work was completed while  the author was visiting  the Vietnam Institute for Advanced Study in Mathematics (VIASM) in the year 2014   for his postdoctoral fellowship. He thanks VIASM for financial support and hospitality.
 
  The work was also supported in part by the grant-aided research project of the Vietnam Ministry of Education and Training.
  
  Finally, the author  would like to express his warmest thanks to the referees for the careful reading
and detailed comments with many helpful suggestions.



\begin{thebibliography}{99}


\bibitem{BH04} Banyaga, A.,   Hurtubise, D.:  Lectures on Morse Homology.  Springer-Verlag, Berlin (2004)
 \bibitem{Bi07} Bivi\`{a}-Ausina, C.: Injectivity of real polynomial maps and Lojasiewicz exponent at infinity. Math. Z. \textbf{257}, 745-767 (2007)
\bibitem{DHP14} Đinh, S.T.,  Hà, H.V.,  Phạm, T.S.:  A Frank-Wolfe type theorem for nondegenerate polynomial programs. Math. Program., Ser. A  \textbf{147}, 519-538  (2014) 
\bibitem{La01} Lasserre, J.-B.: Global optimization with polynomials and the problem of moments. SIAM J. Optim. \textbf{11},   796-817 (2001) 
\bibitem{La10}   Lasserre, J.-B.:   Moments, positive polynomials and their applications.   Imperial College Press Optimization Series, 1. Imperial College Press, London (2010)
 \bibitem{Lau09} Laurent, M.: Sums of squares, moment matrices and optimization over polynomials.  In: Emerging Applications of Algebraic Geometry, IMA Vol. Math. Appl. \textbf{149}, M. Putinar and S. Sullivant, eds., pp. 157-270.  Springer, New York,  (2009)
\bibitem{Le}  L\^{e}, C.-T.:  Characterization of nonnegative polynomials in two variables on compact basic semi-algebraic sets via Newton diagrams. Quy Nhon University Journal of Science VIII (1), 5-16 (2014)
 
\bibitem{Mar03}  Marshall, M.:  Optimization of polynomial functions. Canad. Math. Bull. \textbf{46}, 575-587 (2003) 
\bibitem{Mar08}   Marshall, M.: Positive polynomials and Sums of squares.  Mathematical Surveys and Monographs, \textbf{146},  American Mathematical Society (2008)
\bibitem{Mar09} Marshall, M.: Representation of nonnegative polynomials, degree bounds and applications to optimization.   Canad. J. Math \textbf{61}, 205-221  (2009)
\bibitem{NDP07} Nie, J., Demmel, J.,   Powers, V.:  Representations of positive polynomials on noncompact semialgebraic sets via KKT ideals.  J. Pure Appl. Algebra  \textbf{209} (1), 189–200 (2007)
\bibitem{NDS06} Nie, J., Demmel, J.,  Sturmfels, B.:  Minimizing polynomials via sum of squares over the gradient ideal.  Math. Program., Ser. A \textbf{106}, 587–606 (2006)
\bibitem{NW} Nocedal, J.,  Wright, S.J.:  Numerical optimization. Springer Series in Operations Research,  New York: Springer-Verlag
 (1999)
\bibitem{PS03} Parrilo, P.A.,  Sturmfels, B.: Minimizing polynomial functions.  In: Algorithmic and Quantitative Real Algebraic Geometry, DIMACS Ser. Discrete Math. Theoret. Comput. Sci. \textbf{60}, S. Basu and L. Gonzalez-Vega, eds.,  pp. 83–99. Amer. Math. Soc., Providence, RI (2003)
 \bibitem{Pu} Putinar, M.: Positive polynomials on compact semi-algebraic sets.   Indiana Univ. Math. J. {\bf 42}, no. 3, 969-984 (1993)
\bibitem{Sch}  Schm\"udgen, K.:   The K-moment problem for compact semi-algebraic sets.   Math. Ann. {\bf 289}, 203-206  (1991)
\bibitem{Schei03} Scheiderer,  C.: Sums of squares on real algebraic curves.   Math. Z. {\bf 245}, 725-760 (2003)
 
\bibitem{Schei05} Scheiderer, C.:  Distinguished representations of nonnegative polynomials.   J. Algebra \textbf{289}), no. 2, 558–573  (2005)
  \bibitem{Ste96} Stengle, G.:   Complexity Estimates for the Schm\"{u}dgen Positivstellensatz.  J. Complexity \textbf{12},  167-174 (1996) 
\end{thebibliography}
\end{document}